\documentclass{article}

\usepackage{todonotes}

\usepackage{amsmath,amssymb, amscd, amsthm}
\usepackage{breqn}
\usepackage{mathtools}
\usepackage{autobreak}
\usepackage[pagebackref=true]{hyperref}
\usepackage{listings}
\usepackage{tabularx}

\usepackage{dynkin-diagrams}

\usepackage{breqn}

\usepackage{url}
\hypersetup{breaklinks=true}

\usepackage{breakurl}

\allowdisplaybreaks
\numberwithin{equation}{section}

\theoremstyle{definition}
\newtheorem{definition}{Definition}[section]

\newtheorem{proposition}[definition]{Proposition}

\newtheorem{example}[definition]{Example}

\DeclareMathOperator{\ch}{ch}
\DeclareMathOperator{\td}{td}
\DeclareMathOperator{\Sym}{Sym}

\usepackage{geometry}

\title{A Sage package for calculation of elliptic genera of homogeneous spaces and their complete intersections}
\author{Kobayashi Kenta}
\date{}

\begin{document}
\maketitle

\begin{abstract}
  We present a SageMath package for calculating elliptic genera of homogeneous spaces and their complete intersections.
  This includes the calculation of the basis of weak Jacobi forms,
  Chern numbers of homogeneous spaces and their complete intersections,
  and their elliptic genera.
  In this paper,
  we show how to use this package with concrete examples.
  This package is developed for the calculations in arXiv:2306.10650.
\end{abstract}

\section{Introduction}
For a compact complex manifold $M$,
the two-variable elliptic genus is defined as follows\cite{MR985505, MR1189136, arxiv:math/0405232, MR1048541}:
\begin{align}
  \chi(M, q, y)
   & =
  \int_M
  \ch(\mathcal{ELL}_{q,y})
  \td(M)
  \in
  y^{-d/2} \cdot \mathbb{C}[y^\pm] [[q]]
\end{align}
where
\begin{align}
  \mathcal{ELL}_{q, y}
   & \coloneqq
  y^{- \frac{d}{2}}
  \bigotimes_{n \geq 1}^{\infty}
  \left(
  \wedge_{-yq^{n-1}} T_M^*
  \otimes
  \wedge_{-y^{-1}q^n} T_M
  \otimes
  \Sym_{q^n} T_M^*
  \otimes
  \Sym_{q^n}T_M
  \right),     \\
  \Sym_q V
   & \coloneqq
  \bigoplus_{n \geq 0} ^{\infty}
  \left(
  \Sym^n V
  \right)
  \cdot
  q^n,         \\
  \wedge_q V
   & \coloneqq
  \bigoplus_{n \geq 0}^{\infty}
  \left(
  \wedge^n V
  \right)
  \cdot
  q^n.
\end{align}
If $M$ is a Calabi--Yau $d$-fold,
then the elliptic genus is a weak Jacobi form of weight 0 and index $d/2$ \cite{MR1757003}.
The specialization of the elliptic genus at $q=0$ gives the $\chi_y$ genus.

We developed a SageMath\cite{sagemath} package to calculate such elliptic genera of
homogeneous spaces and their complete intersections.
In this paper,
we explain how to use this.
This package consists of three modules.

\begin{enumerate}
  \item \lstinline|weak_Jacobi_form| provides functions to compute the basis of the space of weak Jacobi forms of even weight.
  \item \lstinline|homogeneous_space| provides functions to compute the cohomology of homogeneous spaces and their complete intersections.
  \item \lstinline|elliptic_genus| provides functions to compute the elliptic genus.
\end{enumerate}

The module \lstinline|weak_Jacobi_form| contains functions to compute the basis of the space of weak Jacobi forms,
and is entirely independent from the other parts.
It contains functions called \lstinline|basis_integral(weight, index)|
(resp. \lstinline|basis_half_integral(weight, index)|),
which return a \lstinline|list| of the basis of the space of weak Jacobi forms with even weight and integral
(resp. half-integral) index.

The module \lstinline|homogeneous_space| contains classes
representing homogeneous spaces and their complete intersections,
and provides methods and global functions to calculate their cohomology.
The included classes are \lstinline|HomogeneousSpace|,
\lstinline|CompletelyReducibleEquivariantVectorBundle|,
\lstinline|IrreducibleEquivariantVectorBundle|,
\lstinline|CompleteIntersection|,
and their abstract superclasses.
These respectively represent homogeneous spaces of quotients by parabolic subgroups, completely reducible equivariant vector bundles over such homogeneous spaces, irreducible equivariant vector bundles,
and complete intersection submanifolds formed by the zeros of general sections of these equivariant vector bundles.
Methods to calculate their
Chern classes,
Chern characters,
and Todd classes,
are also provided.
Moreover,
a global function \lstinline|chern_number(manifold, degrees)| is available,
which calculates the Chern numbers of manifolds.

The module \lstinline|elliptic_genus| contains functions for calculating the elliptic genus.
The function \lstinline|elliptic_genus_chernnum(dim, k)| expresses the terms up to \lstinline|q^k| of the elliptic genus of a \lstinline|dim|-dimensional manifold using Chern numbers.
The function \lstinline|elliptic_genus(manifold, k, option="symbolic")| calculates the terms up to \lstinline|q^k| of the elliptic genus of the object \lstinline|manifold| representing a manifold.

In Section 2,
we explain how to install the package using pip and provide some simple usage examples.
We discuss in more detail how to use the module \lstinline|weak_Jacobi_form| in Section 3,
the module \lstinline|homogeneous_space| in Section 4,
and the module \lstinline|elliptic_genus| in Section 5.

A Jupyter Notebook file containing all examples in this paper is available as an ancillary file to this paper.

This package is developed to prove the coincidence of the elliptic genera of a pair of
complete intersection Calabi--Yau $17$-folds in $F_4$-Grassmannians in \cite{arxiv:2306.10650}.

\section{Getting started}

\subsection{How to install}
The EllipticGenus package,
which is explained in this paper,
is available on GitHub (\burl{https://github.com/topostaro/EllipticGenus}),
and can be installed using the \lstinline|pip| command.

\begin{lstlisting}
  sage --pip install git+https://github.com/topostaro/EllipticGenus.git
\end{lstlisting}

\subsection{The first step}
Here we will compute the elliptic genus of the quintic Calabi--Yau 3-fold inside the $4$-dimensional complex projective space $\mathbb{P}^4$.

First,
we construct the $4$-dimensional complex projective space $\mathbb{P}^4$. The manifold $\mathbb{P}^4$ can be obtained by taking the quotient of the complex Lie group of type $A_4$
by the parabolic subgroup corresponding to the crossed Dynkin diagram $\dynkin[label]{A}{x***}$.
\lstinline|ParabolicSubgroup| is a class representing the parabolic subgroup.
The first argument of the constructor is the Cartan type of the Lie group that contains the parabolic subgroup,
the second argument is the Cartan type of the vertices that are not crossed out \footnote{Although this is determined from the first and third argument,
  at present it needs to be specified by hand.},
and the third argument is the list of the crossed-out vertex
\footnote{In SageMath, the vertices of the Dynkin diagram are numbered from 1 in order.}.
\lstinline|HomogeneousSpace| is a class representing a homogeneous space.
When an object of \lstinline|ParabolicSubgroup| is given as an argument to the constructor,
the homogeneous space obtained by taking the quotient by that parabolic subgroup is returned.

\begin{lstlisting}
  sage: from homogeneous_space import ParabolicSubgroup, HomogeneousSpace
  sage: P = ParabolicSubgroup(CartanType('A4'), CartanType('A3'), [1])
  sage: Proj4 = HomogeneousSpace(P)
\end{lstlisting}

The line bundle $\mathcal{O}(5)$ over $\mathbb{P}^4$ can be constructed as follows. The class \lstinline|IrreducibleEquivariantVectorBundle| represents an irreducible equivariant vector bundle over a homogeneous space.
Its constructor takes as its first argument an object of \lstinline|HomogeneousSpace| as the base space,
and as its second argument the highest weight of the corresponding representation of the parabolic subgroup.
The highest weight is represented as a coefficient vector with the fundamental weights as the basis.

\begin{lstlisting}
  sage: from homogeneous_space import IrreducibleEquivariantVectorBundle
  sage: L = IrreducibleEquivariantVectorBundle(Proj4, (5, 0, 0, 0))
\end{lstlisting}

The quintic Calabi--Yau 3-fold can be obtained as the zero locus of a general section of $\mathcal{O}(5)$.
The class \lstinline|CompleteIntersection| represents complete intersections
that are the zero loci of general sections of equivariant vector bundles
over homogeneous spaces.
Its constructor takes an equivariant vector bundle as an argument.

\begin{lstlisting}
  sage: from homogeneous_space import CompleteIntersection
  sage: Quintic = CompleteIntersection(L)
\end{lstlisting}

The elliptic genus can be computed by calling the function \lstinline|elliptic_genus|.
The first argument of \lstinline|elliptic_genus| is the manifold for which you want to compute the elliptic genus,
and it calculates up to the degree of $q$ given as the second argument.

\begin{lstlisting}
  sage: from elliptic_genus import elliptic_genus
  sage: elliptic_genus(Quintic, 2)
  -100*y - 100*y^2 + (100*y^-1 - 100*y - 100*y^2 + 100*y^4)*q + (100*y^-2 + 100*y^-1 - 200*y - 200*y^2 + 100*y^4 + 100*y^5)*q^2 + O(q^3)
\end{lstlisting}

You can view the documentation for each function and class by typing something like
\newline
\lstinline|elliptic_genus?| in SageMath's interactive mode.

\begin{lstlisting}
  sage: elliptic_genus?
  Signature:
  elliptic_genus(
      manifold: homogeneous_space.interfaces.AlmostComplexManifold,
      k: int,
      option='symbolic',
  )
  Docstring:
    Return the elliptic genus of the argument "manifold" with the terms
    of q variable up to degree "k" multiplied by y^{dim/2} by
    integration with an option "option".

    INPUT:

    * "manifold" -- object of "AlmostComplexManifold" -- this function
      returns the elliptic genus of "manifold"

    * "k" -- integer

    * "option" -- string -- this specifying the integration option.

    OUTPUT:

    the elliptic genus of the argument "manifold" with the terms of q
    variable up to degree "k" multiplied by y^{dim/2}.

    EXAMPLE:
      ...
\end{lstlisting}

\section{Basis of weak Jacobi forms}

Weak Jacobi forms are defined as follows:

\begin{definition}[\cite{MR781735}]
  A weak Jacobi form of weight $k \in \mathbb{Z}$ and index $m \in \frac{1}{2}\mathbb{Z}$ is a function $\tilde{\phi} \left( \tau, z \right)$ with Fourier series expansion
  \begin{equation}
    \tilde{\phi} \left( \tau, z \right)
    =
    \sum_{n \geq 0}
    \sum_{r}
    c\left( n, r \right)
    q^ny^r
    =
    \sum_{n \geq 0}
    \sum_{r}
    c\left( n, r \right)
    \mathrm{e}^{2 \pi \mathrm{i} n \tau}
    \mathrm{e}^{2 \pi \mathrm{i} r z}
  \end{equation}
  which satisfies the following equations:
  \begin{align}
    \tilde{\phi}
    \left(
    \frac{a \tau + b}{c \tau + d},
    \frac{z}         {c \tau + d}
    \right)
     & =
    (c \tau + d)^k
    \mathrm{e}^{\frac{2 \pi \mathrm{i} cz^2}{c \tau + d}}
    \tilde{\phi}
    \left( \tau, z \right),
     & \begin{pmatrix}
         a & b \\
         c & d
       \end{pmatrix}
     & \in \mathrm{SL}_2{\mathbb{Z}}, \\
    \tilde{\phi}
    \left(
    \tau,
    z + \lambda \tau + \mu
    \right)
     & =
    (-1)^{2m(\lambda + \mu)}
    \mathrm{e}
    ^{-2 \pi \mathrm{i} m(\lambda^2 \tau + 2 \lambda z)}
    \tilde{\phi}
    \left( \tau, z \right),
     & \left( \lambda, \mu \right)
     & \in \mathbb{Z}^2.
  \end{align}
  The space of weak Jacobi forms of weight $k$ and index $m$ is denoted by
  $\tilde{\mathrm{J}}_{k, m}$.
\end{definition}

The module \lstinline|weak_Jacobi_form| provides the following two functions.

\begin{enumerate}
  \item
        \lstinline|basis_integral(weight, index)| returns the basis of the space spanned by weak Jacobi forms with even weight and integral index.
  \item
        \lstinline|basis_half_integral(weight, double_index)| returns the basis of the space spanned by weak Jacobi forms with even weight and half-integral index \lstinline|double_index/2|.
\end{enumerate}

\subsection{\lstinline|basis_integral(weight, index)|}

Proposition \ref{pr:basis} below shows that
the space of weak Jacobi forms with even weight and integral index
is a polynomial ring.

\begin{proposition}[{\cite[Theorem 9.3]{MR781735}}] \label{pr:basis}
  The space $\tilde{J}_{\text{even}, *}$ of weak Jacobi forms with even weight and integral index
  is isomorphic to $\mathrm{M}_*[\tilde{\phi}_{-2, 1}, \tilde{\phi}_{0, 1}]$.
  Here,
  $\mathrm{M}_* = \mathbb{C}[\mathrm{E}_4, \mathrm{E}_6]$ is the ring of modular forms
  ($\mathrm{E}_4$ and $\mathrm{E}_6$ are the normalized Eisenstein series of weight $4$ and $6$),
  and
  \begin{align}
    \tilde{\phi}_{0, 1}(q, y)
     & = 4 \left\{\left(\frac{\theta_2(q, y)}{\theta_2(q, 1)}\right)^2
    + \left(\frac{\theta_3(q, y)}{\theta_3(q, 1)}\right)^2
    + \left(\frac{\theta_4(q, y)}{\theta_4(q, 1)}\right)^2\right\},    \\
    \tilde{\phi}_{-2, 1}(q, y)
     & = - \frac{\theta_1(q, y)^2}{\eta(q)^6}.
  \end{align}
  The function $\tilde{\phi}_{0, 1}(q, y)$ (resp. $\tilde{\phi}_{-2, 1}(q, y)$) is a weak Jacobi form of weight $0$ and index $1$ (resp. weight $-2$ and index $1$).
\end{proposition}

The function \lstinline|basis_integral(weight, index)| lists
products of the generators
with the correct weight and index.

\begin{example}
  The basis of the space of weak Jacobi forms of weight $0$ and index $3$ can be obtained as follows.

  \begin{lstlisting}
    sage: from weak_Jacobi_form import basis_integral
    sage: basis_integral(0, 3)
    [(y^-3 + 30*y^-2 + 303*y^-1 + 1060 + 303*y + 30*y^2 + y^3) + ... + O(q^7),
     (y^-3 + 6*y^-2 - 33*y^-1 + 52 - 33*y + 6*y^2 + y^3) +  ... + O(q^7),
     (y^-3 - 6*y^-2 + 15*y^-1 - 20 + 15*y - 6*y^2 + y^3) + ... + O(q^7)]
  \end{lstlisting}

  The output represents the set
  \begin{equation}
    \left\{
    \tilde{\phi}_{0, 1}^3,~
    \tilde{\phi}_{0, 1}\tilde{\phi}_{-2, 1}^2\mathrm{E}_4,~
    \tilde{\phi}_{-2, 1}^3 \mathrm{E}_6
    \right\}.
  \end{equation}
\end{example}

\subsection{\lstinline|basis_half_integral(weight, double_index)|}

This function \lstinline|basis_half_integral(weight, double_index)| computes the basis of the space spanned by weak Jacobi forms with even weight and half-integral index.
The computation is based on Proposition \ref{pr:half-integral} below.

\begin{proposition}[\cite{MR1757003}] \label{pr:half-integral}
  For any positive odd integer $d$,
  there exists an isomorphism
  \begin{equation}
    \tilde{J}_{0, \frac{(d-3)}{2}} \rightarrow \tilde{J}_{0, \frac{d}{2}}, \qquad
    \phi \mapsto \tilde{\phi}_{0, \frac{3}{2}}(q, y) \phi
  \end{equation}
  where
  \begin{equation}
    \tilde{\phi}_{0, \frac{3}{2}}(q, y) \coloneqq \frac{\theta_1(q, y^2)}{\theta_1(q, y)}.
  \end{equation}
\end{proposition}

When $d$ is even,
it returns the same result as \lstinline|basis_integral|.
When $d$ is odd,
for the convenience of representation in SageMath,
it returns the value of the basis of the weak Jacobi form space multiplied by $y^{1/2}$.

\begin{example}
  The space of weak Jacobi forms of weight $2$ and index $5/2$ can be obtained as follows.

  \begin{lstlisting}
    sage: from weak_Jacobi_form import basis_half_integral
    sage: basis_half_integral(2, 5)
    [(y^-1 - 1 - y + y^2) + (-y^-3 + 246*y^-1 - 245 - 245*y + 246*y^2 - y^4)*q + ... + O(q^7)]
  \end{lstlisting}

\end{example}

\section{Chern classes of homogeneous spaces}

The module \lstinline|homogeneous_space| provides classes and functions to compute the cohomology of homogeneous spaces formed by the quotient by a parabolic subgroup and complete intersections of the zero locus of a general section of a covariant vector bundle over them.

\subsection{\lstinline|ParabolicSubgroup|}

A subgroup $P$ of a complex reductive Lie group $G$ is \emph{parabolic} if $G/P$ is projective.
The set of conjugacy classes of parabolic subgroups is bijective
with the set of subsets of the set of simple roots of $G$.
Through this correspondence,
a (conjugacy class of) parabolic subgroup can be represented
by a \emph{crossed Dynkin diagram};
nodes corresponding to simple roots not contained in the subset
are crossed out.
In this package,
this is represented by the class \lstinline|ParabolicSubgroup|.
The first argument of its constructor is the Cartan type of the ambient Lie group $G$,
the second argument is the Cartan type formed by the non-crossed-out nodes,
and the third argument is the list consisting of crossed-out nodes.%
\footnote{As mentioned in the footnote of Section 2.2,
  the second argument is redundant,
  providing an opportunity for improvement.}
When all nodes are crossed out,
put \lstinline|None| to the second argument.

\begin{example}
  The parabolic subgroup of the simple Lie group of type $F_4$ corresponding to the crossed Dynkin diagram $\dynkin[label]{F}{*x**}$ can be obtained as follows.
  \vspace{0.2cm}
  \begin{lstlisting}
    sage: from homogeneous_space import ParabolicSubgroup
    sage: P = ParabolicSubgroup(CartanType('F4'), CartanType('A1xA2'), [2])
    sage: P.dynkin_diagram()
    O---X=>=O---O
    1   2   3   4
    F4 with node 2 marked
  \end{lstlisting}
\end{example}

\begin{example}
  The parabolic subgroup of the simple Lie group of type $A_2$ corresponding to the crossed Dynkin diagram $\dynkin[label]{A}{xx}$ can be obtained as follows.
  \vspace{0.2cm}
  \begin{lstlisting}
    sage: from homogeneous_space import ParabolicSubgroup
    sage: P = ParabolicSubgroup(CartanType('A2'), None, [1, 2])
    sage: P.dynkin_diagram()
    X---X
    1   2
    A2 with nodes (1, 2) marked
  \end{lstlisting}
\end{example}

The set of simple roots of the parabolic subgroup $P$
is the subset of the set of simple roots of $G$.
Positive roots of $P$ are roots of $G$
which are positive linear combinations of simple roots of $P$.

\begin{example}

  Simple roots and positive roots of the parabolic subgroup $P$ of type $A_4$ corresponding to the crossed Dynkin diagram \dynkin[label]{A}{**x*} can be computed as follows.
  In the output of \lstinline|P.simple_roots()|, the simple root \lstinline|(0, 0, 1, -1, 0)| corresponding to the third vertex is omitted from \lstinline|CartanType('A4').simple_roots()|.
  In the output of \lstinline|P.positive_roots()|, any roots that is still a root after subtracting \lstinline|(0, 0, 1, -1, 0)| is omitted from \lstinline|CartanType('A4').positive_roots()|.

  \begin{lstlisting}
    sage: from homogeneous_space import ParabolicSubgroup
    sage: P = ParabolicSubgroup(CartanType('A4'), CartanType('A2xA1'), [3])
    sage: P.dynkin_diagram()
    O---O---X---O
    1   2   3   4
    A4 with node 3 marked
    sage: P.simple_roots()
    [(1, -1, 0, 0, 0), (0, 1, -1, 0, 0), (0, 0, 0, 1, -1)]
    sage: P.positive_roots()
    [(1, -1, 0, 0, 0), (0, 1, -1, 0, 0), (0, 0, 0, 1, -1), (1, 0, -1, 0, 0)]
  \end{lstlisting}

\end{example}

The method \lstinline|weight_multiplicities(weight)|
takes a weight
(presented as a vector with the set of fundamental weights as a basis)
as an argument,
and returns the set of weights of the highest weight representation of the parabolic subgroup
as a \lstinline|dict| object;
a key represents a weight, and the value represents the multiplicity of that key.%
\footnote{
  The calculation is valid
  only if the argument is integral and $\mathfrak{p}$-dominant.
}

\begin{example}
  The weights of the highest weight representations
  with highest weights $\omega_1 + 3\omega_3 + \omega_4$ and
  $\omega_2 - \omega_3$
  of the parabolic subgroup $P$
  corresponding to the crossed Dynkin diagram \dynkin[label]{A}{**x*}
  is calculated as follows.

  \begin{lstlisting}
    sage: from homogeneous_space import ParabolicSubgroup
    sage: P = ParabolicSubgroup(CartanType('A4'), CartanType('A2xA1'), [3])
    sage: P.weight_multiplicities((1, 0, 3, 1))
    {(5, 4, 4, 1, 0): 1,
     (4, 5, 4, 1, 0): 1,
     (5, 4, 4, 0, 1): 1,
     (4, 4, 5, 1, 0): 1,
     (4, 5, 4, 0, 1): 1,
     (4, 4, 5, 0, 1): 1}
    sage: P.weight_multiplicities((0, 1, -1, 0))
    {(0, 0, -1, 0, 0): 1, (0, -1, 0, 0, 0): 1, (-1, 0, 0, 0, 0): 1}
  \end{lstlisting}

\end{example}

\subsection{\lstinline|HomogeneousSpace|}

In this package, the class \lstinline|HomogeneousSpace| represents the complex projective manifold obtained by taking the quotient by a parabolic subgroup.

\begin{example}

  If $G \coloneqq \mathrm{SL}(5, \mathbb{C})$ is a Lie group of type $A_4$
  and $P$ is the parabolic subgroup corresponding to the crossed Dynkin diagram \dynkin[label]{A}{**x*},
  then $G/P$ is the complex Grassmannian $\mathrm{Gr}(3, 5)$.

  \begin{lstlisting}
    sage: from homogeneous_space import ParabolicSubgroup, HomogeneousSpace
    sage: P = ParabolicSubgroup(CartanType('A4'), CartanType('A2xA1'), [3])
    sage: Gr35 = HomogeneousSpace(P)
  \end{lstlisting}

\end{example}

The following is a list of methods for the class \lstinline|HomogeneousSpace|.

\vspace{0.4cm}

\noindent
\begin{tabularx}{\textwidth}{p{4.5cm}lX}
  \hline
  Method                                              & Return Type              & Description                                                                                                                                                                                                                                 \\
  \hline
  \lstinline|dimension(self)|                         & \lstinline|int|          & Returns the dimension.                                                                                                                                                                                                                      \\
  \hline
  \lstinline|tangent_bundle(self)|                    & \lstinline|VectorBundle| & Returns the tangent bundle.                                                                                                                                                                                                                 \\
  \hline
  \lstinline|cotangent_bundle(self)|                  & \lstinline|VectorBundle| & Returns the cotangent bundle.                                                                                                                                                                                                               \\
  \hline
  \lstinline|chern_classes(self)|                     & \lstinline|list|         & Returns the Chern classes as a list by degree.                                                                                                                                                                                              \\
  \hline
  \lstinline|todd_classes(self)|                      & \lstinline|list|         & Returns the Todd classes as a list by degree.                                                                                                                                                                                               \\
  \hline
  \lstinline|integration(self, f, option="symbolic")| & \lstinline|int|          & Performs integration of the cohomology class \lstinline|f| using the method specified by \lstinline|option|.                                                                                                                                \\
  \hline
  \lstinline|symbolic_integration_by_|
  \lstinline|localization(self, f)|                   & \lstinline|int|          & Calculates the integral of the cohomology class \lstinline|f| by a symbolic computation using the localization formula in equivariant cohomology. It is called if \lstinline|option| for \lstinline|integration| is \lstinline|"symbolic"|. \\
  \hline
  \lstinline|numerical_integration_by_|
  \lstinline|localization(self, f)|                   & \lstinline|int|          & Calculates the integral by a numerical computation using the localization formula. It is called if \lstinline|option| for \lstinline|integration| is \lstinline|"numerical"|.                                                               \\
  \hline
\end{tabularx}

\vspace{0.4cm}

To be more precise,
the return class of \lstinline|tangent_bundle|
and \lstinline|cotangent_bundle|
is a subclass of the abstract class
\lstinline|VectorBundle|.
The abstract class
\lstinline|VectorBundle|
is designed to keep track of the information of Chern classes.

The method \lstinline|chern_classes(self)| returns Chern classes.
Chern classes are represented as elements of $\operatorname{rank} G$-variable polynomials.
Let $\Lambda$ be the weight lattice of $G$,
and $V$ be its default ambient space in SageMath.
Then we have
\begin{equation}
  \mathrm{H}^*_G(G/P) \subseteq H^*_G(G/B) \cong \Sym \Lambda \subseteq \Sym V,
\end{equation}
and Chern classes are represented as elements of $\Sym V$.
$\Sym V$ has generators \lstinline|x0, ..., xn| corresponding to the standard basis of $V$.

The method \lstinline|chern_classes(self)| returns Chern classes
as a list of the form \lstinline|[c0, c1, ... , cn]|,
where
\lstinline|ci| is the $i$-th Chern class
and \lstinline|n|
is the dimension of the manifold.

\begin{example}
  The Chern classes of $\mathrm{Gr}(3, 5)$ and the integral of the top Chern class is calculated as follows.

  \begin{lstlisting}
    sage: from homogeneous_space import ParabolicSubgroup, HomogeneousSpace
    sage: P = ParabolicSubgroup(CartanType('A4'), CartanType('A2xA1'), [3])
    sage: Gr35 = HomogeneousSpace(P)
    sage: Gr35.dimension()
    6
    sage: Gr35.chern_classes()
    [1,
     2*x0 + 2*x1 + 2*x2 - 3*x3 - 3*x4,
     x0^2 + 4*x0*x1 + x1^2 + 4*x0*x2 + 4*x1*x2 + x2^2 - 5*x0*x3 - 5*x1*x3 - 5*x2*x3 + 3*x3^2 - 5*x0*x4 - 5*x1*x4 - 5*x2*x4 + 9*x3*x4 + 3*x4^2,
     ... ]
    sage: Gr35.integration(Gr35.chern_classes()[Gr35.dimension()])
    10
  \end{lstlisting}

\end{example}

The method \lstinline|integration(self, f, option="symbolic")|
provides two ways
to calculate the integral
specified by \lstinline|option|.
Localization formula describes the integral
as the sum of rational functions,
which simplifies to an integer.
If \lstinline|option| is \lstinline|"symbolic"|
(which is the default value),
the calculation is performed symbolically,
and always returns the correct answer.
However,
the calculation is time-consuming if the dimension of the manifold is high.
If \lstinline|option| is \lstinline|"numerical"|,
the integral is calculated by substituting random real values into each variable
and performing numerical calculations.
While there is a very small chance of giving a wrong answer,
the calculation is much faster in general.

\begin{example}
  The integral of the $6$-th power of the first Chern class $\mathrm{c}_1$ of the Grassmannian $\operatorname{Gr}(3, 5)$ is computed using both options below.

  \begin{lstlisting}
    sage: from homogeneous_space import ParabolicSubgroup, HomogeneousSpace
    sage: P = ParabolicSubgroup(CartanType('A4'), CartanType('A2xA1'), [3])
    sage: Gr35 = HomogeneousSpace(P)
    sage: Gr35.integration(Gr35.chern_classes()[1]^6, option="symbolic")
    78125
    sage: Gr35.integration(Gr35.chern_classes()[1]^6, option="numerical")
    78125
  \end{lstlisting}

\end{example}

The tangent bundle and the cotangent bundle can be computed using the methods \lstinline|tangent_bundle(self)| and \lstinline|cotangent_bundle(self)|,
respectively.
The base class for vector bundles,
\lstinline|VectorBundle|,
provides methods for computing the Chern classes.
The method \lstinline|chern_classes(self)| of \lstinline|HomogeneousSpace| computes by calling \lstinline|tangent_bundle(self).chern_classes()|.

\begin{example}
  In this example,
  we calculate the Chern classes of the tangent and cotangent bundles of the $4$-dimensional projective space $\mathbb{P}^4$.

  \begin{lstlisting}
    sage: from homogeneous_space import *
    sage: P = ParabolicSubgroup(CartanType('A4'), CartanType('A3'), [1])
    sage: Proj4 = HomogeneousSpace(P)
    sage: Proj4.tangent_bundle().chern_classes()
    [1,
     4*x0 - x1 - x2 - x3 - x4,
     6*x0^2 - 3*x0*x1 - 3*x0*x2 + x1*x2 - 3*x0*x3 + x1*x3 + x2*x3 - 3*x0*x4 + x1*x4 + x2*x4 + x3*x4,
     ...]
    sage: Proj4.cotangent_bundle().chern_classes()
    [1,
     -4*x0 + x1 + x2 + x3 + x4,
     6*x0^2 - 3*x0*x1 - 3*x0*x2 + x1*x2 - 3*x0*x3 + x1*x3 + x2*x3 - 3*x0*x4 + x1*x4 + x2*x4 + x3*x4,
     ...]
  \end{lstlisting}

\end{example}

\subsection{\lstinline|CompletelyReducibleEquivariantVectorBundle|}

The class \lstinline|CompletelyReducibleEquivariantVectorBundle| is a class representing a completely reducible equivariant vector bundle over a homogeneous space.
The category of G-equivariant vector bundles over the homogeneous space $G/P$
and the category of representations of $P$ are equivalent.
This equivalence sends a representation $V$ of $P$
to the $G$-equivariant vector bundle
$
  \mathcal{E}_V \coloneqq G \times_P V.
$

Since irreducible representations of $P$ correspond bijectively
to irreducible representations of the Levi part of $P$,
a completely reducible representation can be associated
with a sequence of integral $\mathfrak{p}$-dominant weights
by taking the highest weights of its irreducible components.

The constructor of this class takes an object of \lstinline|HomogeneousSpace|,
which serves as the base space,
as the first argument.
As the second argument,
it takes a \lstinline|list| of highest weights of the irreducible components contained in $V$.
It returns an object representing the corresponding equivariant vector bundle $\mathcal{E}_V$.

\begin{example}

  The Chern classes of the rank-$2$ vector bundle $\mathcal{O}(2) \oplus \mathcal{O}(3)$ on the $4$-dimensional projective space $\mathbb{P}^4$ can be computed as follows.

  \begin{lstlisting}
    sage: from homogeneous_space import *
    sage: P = ParabolicSubgroup(CartanType('A4'), CartanType('A3'), [1])
    sage: Proj4 = HomogeneousSpace(P)
    sage: E = CompletelyReducibleEquivariantVectorBundle(Proj4, [(2, 0, 0, 0), (3, 0, 0, 0)])
    sage: E.chern_classes()
    [1, 5*x0, 6*x0^2, 0, 0]
  \end{lstlisting}

\end{example}

The following is a list of methods for the class \lstinline|CompletelyReducibleEquivariantVectorBundle|.

\vspace{0.4cm}

\noindent
\begin{tabularx}{\textwidth}{p{4.5cm}lX}
  Method                            & Return Type                  & Description                                                                                                                                                                                                             \\
  \hline
  \lstinline|rank(self)|            & \lstinline|int|              & Returns the rank.                                                                                                                                                                                                       \\
  \hline
  \lstinline|base(self)|            & \lstinline|HomogeneousSpace| & Returns the base space.                                                                                                                                                                                                 \\
  \hline
  \lstinline|is_irreducible(self)|  & \lstinline|bool|             & Returns whether it is irreducible.                                                                                                                                                                                      \\
  \hline
  \lstinline|chern_classes(self)|   & \lstinline|list|             & Returns the Chern classes as a list by degree.                                                                                                                                                                          \\
  \hline
  \lstinline|chern_character(self)| & \lstinline|list|             & Returns the Chern character as a list by degree.                                                                                                                                                                        \\
  \hline
  \lstinline|todd_classes(self)|    & \lstinline|list|             & Returns the Todd classes as a list by degree.                                                                                                                                                                           \\
  \hline
  \lstinline|dual(self)|            & \lstinline|VectorBundle|     & Returns the dual bundle.                                                                                                                                                                                                \\
  \hline
  \lstinline|__add__(self, other)|  & \lstinline|VectorBundle|     & Overloads the \lstinline|+|-operator and returns the direct sum of \lstinline|self| and \lstinline|other|. There also exists a global function version, \lstinline|direct_sum(vector_bundle1, vector_bundle2)|.         \\
  \hline
  \lstinline|__mul__(self, other)|  & \lstinline|VectorBundle|     & Overloads the \lstinline|*|-operator and returns the tensor product of \lstinline|self| and \lstinline|other|. There also exists a global function version, \lstinline|tensor_product(vector_bundle1, vector_bundle2)|. \\
  \hline
\end{tabularx}

\vspace{0.4cm}

Recall that
the Chern character of the vector bundle $\mathcal{E}$ of rank $r$ is defined using its Chern roots $x_1, \ldots, x_r$ as follows:
\begin{align}
  \ch \mathcal{E}
   & =
  \sum_{i=1}^{r}
  e^{x_i} \\
   & =
  \sum_{i=1}^{r}
  (1 + x_i + \frac{x_i^2}{2} + \frac{x_i^3}{6} + \frac{x_i^4}{24} + \cdots).
\end{align}
Similarly, the Todd classes are defined as follows:
\begin{align}
  \td \mathcal{E}
   & =
  \prod_{i = 1}^r
  \frac{x_i}{1 - e^{x_i}} \\
   & =
  \prod_{i = 1}^r
  (1 + \frac{x_i}{2} + \frac{x_i^2}{12} - \frac{x_i^4}{720} + \cdots).
\end{align}

\begin{example}
  The Chern character and Todd classes of the line bundle $\mathcal{O}(1)$ on the $4$-dimensional projective space $\mathbb{P}^4$ is computed as follows.
  \begin{lstlisting}
    sage: from homogeneous_space import *
    sage: P = ParabolicSubgroup(CartanType('A4'), CartanType('A3'), [1])
    sage: Proj4 = HomogeneousSpace(P)
    sage: E = CompletelyReducibleEquivariantVectorBundle(Proj4, [(1, 0, 0, 0)])
    sage: E.chern_character()
    [1, x0, 1/2*x0^2, 1/6*x0^3, 1/24*x0^4]
    sage: E.todd_classes()
    [1, 1/2*x0, 1/12*x0^2, 0, -1/720*x0^4]
  \end{lstlisting}

\end{example}

\begin{example}

  The direct sum and the tensor product of the line bundles $\mathcal{O}(2)$ and $\mathcal{O}(3)$ on $\mathbb{P}^4$ is computed as follows.

  \begin{lstlisting}
    sage: from homogeneous_space import *
    sage: P = ParabolicSubgroup(CartanType('A4'), CartanType('A3'), [1])
    sage: Proj4 = HomogeneousSpace(P)
    sage: E1 = CompletelyReducibleEquivariantVectorBundle(Proj4, [(2, 0, 0, 0)])
    sage: E2 = CompletelyReducibleEquivariantVectorBundle(Proj4, [(3, 0, 0, 0)])
    sage: E1.chern_classes()
    [1, 2*x0, 0, 0, 0]
    sage: E2.chern_classes()
    [1, 3*x0, 0, 0, 0]
    sage: (E1 + E2).chern_classes()
    [1, 5*x0, 6*x0^2, 0, 0]
    sage: (E1 * E2).chern_classes()
    [1, 5*x0, 0, 0, 0]
  \end{lstlisting}

\end{example}

There also exist global functions to calculate the Chern classes of the symmetric and anti-symmetric powers of vector bundles.%
\footnote{The functions for calculating the Chern classes of direct sums, tensor products, symmetric and anti-symmetric powers are wrappers of the library \lstinline|chern.lib| \cite{chernlib} included in Singular.}
The function \lstinline|symmetric_power(vector_bundle, k)| is a function to calculate the symmetric power of a vector bundle,
which takes a vector bundle as the first argument and returns the symmetric power specified by the integer \lstinline|k| taken as the second argument.
The function \lstinline|wedge_power(vector_bundle, k)| is a function to calculate the anti-symmetric power of a vector bundle,
which takes a vector bundle as the first argument and returns the anti-symmetric power specified by the integer \lstinline|k| taken as the second argument.

\begin{example}

  The third symmetric power $\Sym^3 \mathcal{O}(2)$ and the second anti-symmetric power $\wedge^2 \mathcal{O}(2)$ of the line bundle $\mathcal{O}(2)$ on $\mathbb{P}^4$ is calculated as follows.

  \begin{lstlisting}
    sage: from homogeneous_space import *
    sage: P = ParabolicSubgroup(CartanType('A4'), CartanType('A3'), [1])
    sage: Proj4 = HomogeneousSpace(P)
    sage: L = CompletelyReducibleEquivariantVectorBundle(Proj4, [(2, 0, 0, 0)])
    sage: symmetric_power(L, 3).chern_classes()
    [1, 6*x0, 0, 0, 0]
    sage: wedge_power(L, 2).chern_classes()
    [1, 0, 0, 0, 0]
  \end{lstlisting}

\end{example}

There also exists a subclass \lstinline|IrreducibleEquivariantVectorBundle| representing irreducible equivariant vector bundles.
The constructor takes a homogeneous space as the first argument,
and the highest weight of an irreducible representation as the second argument.

\begin{example}

  Chern classes of an irreducible equivariant vector bundle of rank $3$ on $\mathbb{P}^4$ is computed as follows.

  \begin{lstlisting}
    sage: from homogeneous_space import *
    sage: P = ParabolicSubgroup(CartanType('A4'), CartanType('A3'), [1])
    sage: Proj4 = HomogeneousSpace(P)
    sage: E = IrreducibleEquivariantVectorBundle(Proj4, (0, 1, 0, 0))
    sage: E.chern_classes()
    [1,
     4*x0 + x1 + x2 + x3 + x4,
     6*x0^2 + 3*x0*x1 + 3*x0*x2 + x1*x2 + 3*x0*x3 + x1*x3 + x2*x3 + 3*x0*x4 + x1*x4 + x2*x4 + x3*x4,
     4*x0^3 + 3*x0^2*x1 + 3*x0^2*x2 + 2*x0*x1*x2 + 3*x0^2*x3 + 2*x0*x1*x3 + 2*x0*x2*x3 + x1*x2*x3 + 3*x0^2*x4 + 2*x0*x1*x4 + 2*x0*x2*x4 + x1*x2*x4 + 2*x0*x3*x4 + x1*x3*x4 + x2*x3*x4,
     x0^4 + x0^3*x1 + x0^3*x2 + x0^2*x1*x2 + x0^3*x3 + x0^2*x1*x3 + x0^2*x2*x3 + x0*x1*x2*x3 + x0^3*x4 + x0^2*x1*x4 + x0^2*x2*x4 + x0*x1*x2*x4 + x0^2*x3*x4 + x0*x1*x3*x4 + x0*x2*x3*x4 + x1*x2*x3*x4]
  \end{lstlisting}

\end{example}

\subsection{\lstinline|CompleteIntersection|}

The class \lstinline|CompleteIntersection| represents
a submanifold obtained as the zero locus
of a general section of a \lstinline|CompletelyReducibleEquivariantVectorBundle|.
The constructor takes a vector bundle as an argument.

\begin{example}

  The zero locus of a general section of the line bundle $\mathcal{O}(4)$ on the $3$-dimensional projective space $\mathbb{P}^3$ gives a quartic K3 surface.

  \begin{lstlisting}
    sage: from homogeneous_space import *
    sage: P = ParabolicSubgroup(CartanType('A3'), CartanType('A2'), [1])
    sage: Proj3 = HomogeneousSpace(P)
    sage: E = IrreducibleEquivariantVectorBundle(Proj3, (4, 0, 0))
    sage: K3 = CompleteIntersection(E)
  \end{lstlisting}

\end{example}

The interface of this class is almost the same as that of \lstinline|HomogeneousSpace|.

\vspace{0.4cm}

\noindent
\begin{tabularx}{\textwidth}{p{4.5cm}lX}
  \hline
  Method                                              & Return Type              & Description                                                 \\
  \hline
  \lstinline|dimension(self)|                         & \lstinline|int|          & Returns the dimension.                                      \\
  \hline
  \lstinline|tangent_bundle(self)|                    & \lstinline|VectorBundle| & Returns the tangent bundle.                                 \\
  \hline
  \lstinline|cotangent_bundle(self)|                  & \lstinline|VectorBundle| & Returns the cotangent bundle.                               \\
  \hline
  \lstinline|chern_classes(self)|                     & \lstinline|list|         & Returns the Chern classes as a list by degree.              \\
  \hline
  \lstinline|todd_classes(self)|                      & \lstinline|list|         & Returns the Todd classes as a list by degree.               \\
  \hline
  \lstinline|integration(self, f, option="symbolic")| & \lstinline|int|          & Performs integration of the cohomology class \lstinline|f|.
  The option specifies the options when calling \lstinline|integration| of \lstinline|HomogeneousSpace| during the calculation process.        \\
  \hline
\end{tabularx}

\vspace{0.4cm}

\begin{example}
  The integral of the second Chern class of a quartic K3 surface is computed as follows.

  \begin{lstlisting}
    sage: from homogeneous_space import *
    sage: P = ParabolicSubgroup(CartanType('A3'), CartanType('A2'), [1])
    sage: Proj3 = HomogeneousSpace(P)
    sage: E = IrreducibleEquivariantVectorBundle(Proj3, (4, 0, 0))
    sage: K3 = CompleteIntersection(E)
    sage: K3.integration(K3.chern_classes()[2], option="numerical")
    24
  \end{lstlisting}

\end{example}

\subsection{\lstinline|chern_number(manifold, degrees, option="symbolic")|}

The global function \lstinline|chern_number(manifold, degrees, option="symbolic")|
can be used to calculate the Chern numbers of \lstinline|HomogeneousSpace|
or \lstinline|CompleteIntersection|.
The first argument should be an object of \lstinline|HomogeneousSpace|
or \lstinline|CompleteIntersection| for which the Chern numbers are to be calculated.
The second argument is a list of degrees for which the Chern numbers are calculated.
The third argument is an option for integration,
which specifies the options when calling the method \lstinline|integration| of the object \lstinline|manifold| which is the first argument.

\begin{example}
  In this example, we calculate the Chern numbers of the complex Grassmannian $\mathrm{Gr}(3, 5)$.
  The following integrals are computed:
  $\int_{\mathrm{Gr}(3, 5)} c_1^6$,
  $\int_{\mathrm{Gr}(3, 5)} c_1 c_2 c_3$,
  $\int_{\mathrm{Gr}(3, 5)} c_6$,
  and
  $\int_{\mathrm{Gr}(3, 5)} c_3 c_4$.
  If the elements of the list do not add up to the dimension of the manifold,
  then the integral is $0$.

  \begin{lstlisting}
    sage: from homogeneous_space import *
    sage: P = ParabolicSubgroup(CartanType('A4'), CartanType('A2xA1'), [3])
    sage: Gr35 = HomogeneousSpace(P)
    sage: chern_number(Gr35, [1] * 6)
    78125
    sage: chern_number(Gr35, [1, 2, 3])
    4275
    sage: chern_number(Gr35, [6])
    10
    sage: chern_number(Gr35, [3, 4])
    0
  \end{lstlisting}

\end{example}

\section{Elliptic genera}

The module \lstinline|elliptic_genus| provides the calculation of the elliptic genus in two stages.

\begin{enumerate}
  \item \lstinline|elliptic_genus_chernnum(dim, k)| expresses the elliptic genus of a \lstinline|dim|-dimensional manifold up to the \lstinline|q^k|-term using Chern numbers.
  \item \lstinline|elliptic_genus(manifold, k, option="symbolic")| calculates the elliptic genus up to the \lstinline|q^k|-term using the Chern number values of the manifold object \lstinline|manifold| with integration specified by \lstinline|option|.
\end{enumerate}

\subsection{\lstinline|elliptic_genus_chernnum(dim, k)|}
The elliptic genus can be expressed in terms of the Chern roots $x_1, \ldots, x_d$ of the tangent bundle as follows:
\begin{align}
  \chi(M, q, y)
   & =
  y^{- \frac{d}{2}}
  \int_M
  \ch\left(
  \bigotimes_{n \geq 0}^{\infty}
  \left(
    \wedge_{-yq^{n-1}} T_M^*
    \otimes
    \wedge_{-y^{-1}q^n} T_M
    \otimes
    \Sym_{q^n} T_M^*
    \otimes
    \Sym_{q^n}T_M
    \right)
  \right)
  \td\left( M \right) \\
   & =
  y^{- \frac{d}{2}}
  \int_M
  \prod_{i = 1}^{d} \prod_{n=1}^\infty
  \frac{(1 - yq^{n-1}e^{-x_i})(1 - y^{-1}q^ne^{x_i})}{(1 - q^ne^{-x_i})(1 - q^ne^{x_i})}
  \cdot
  \frac{x_i}{1 - e^{-x_i}}.
\end{align}
Since the integrand is symmetric in $x_i$,
it can be expressed in terms of the elementary symmetric polynomials of $x_i$. The elementary symmetric polynomials of the Chern roots are the Chern classes of the same degree.
After integration, one obtains a presentation of the elliptic genus using Chern numbers. For the sake of presentation as SageMath objects, \lstinline|elliptic_genus_chernnum(dim, k)| expresses the integrand multiplied by $y^{d/2}$ in terms of Chern classes.

\begin{example}
  The elliptic genus of a 3-dimensional manifold up to first order in $q$ is expressed in terms of Chern numbers as follows.

  \begin{lstlisting}
    sage: from elliptic_genus import elliptic_genus_chernnum
    sage: elliptic_genus_chernnum(3,1)
    1/24*c1*c2 + (-1/24*c1*c2 + 1/2*c3)*y + (-1/24*c1*c2 + 1/2*c3)*y^2 + 1/24*c1*c2*y^3 + ((-1/2*c1^3 + 19/24*c1*c2 - 1/2*c3)*y^-1 + (3/2*c1^3 - 27/8*c1*c2) + (-c1^3 + 31/12*c1*c2 + 1/2*c3)*y + (-c1^3 + 31/12*c1*c2 + 1/2*c3)*y^2 + (3/2*c1^3 - 27/8*c1*c2)*y^3 + (-1/2*c1^3 + 19/24*c1*c2 - 1/2*c3)*y^4)*q + O(q^2)
  \end{lstlisting}

\end{example}

\subsection{\lstinline|elliptic_genus(manifold, k, option="symbolic")|}

The function \lstinline|elliptic_genus(manifold, k, option="symbolic")| computes the elliptic genus by substituting the actual Chern numbers of \lstinline|manifold| into the expression calculated by \lstinline|elliptic_genus_chernnum(dim, k)|.
Due to the inner presentation as an object in SageMath,
\lstinline|elliptic_genus(manifold, k, option="symbolic")| outputs the actual elliptic genus multiplied by $y^{d/2}$, where $d$ is the (complex) dimension of \lstinline|manifold|.
The third argument \lstinline|option| specifies the integration options when calculating the Chern numbers of \lstinline|manifold|.

\begin{example}
  The elliptic genus of a quartic K3 surface multiplied by $y$ is computed as follows.

  \begin{lstlisting}
    sage: from homogeneous_space import *
    sage: from elliptic_genus import elliptic_genus
    sage: P = ParabolicSubgroup(CartanType('A3'), CartanType('A2'), [1])
    sage: Proj3 = HomogeneousSpace(P)
    sage: E = IrreducibleEquivariantVectorBundle(Proj3, (4, 0, 0))
    sage: K3 = CompleteIntersection(E)
    sage: elliptic_genus(K3, 3)
    2 + 20*y + 2*y^2 + (20*y^-1 - 128 + 216*y - 128*y^2 + 20*y^3)*q + (2*y^-2 + 216*y^-1 - 1026 + 1616*y - 1026*y^2 + 216*y^3 + 2*y^4)*q^2 + (-128*y^-2 + 1616*y^-1 - 5504 + 8032*y - 5504*y^2 + 1616*y^3 - 128*y^4)*q^3 + O(q^4)
  \end{lstlisting}

\end{example}

\begin{example}
  The elliptic genus of a Calabi--Yau 3-fold in the $G_2$-flag variety
  from \cite[Table 1.2, No.~5]{arxiv:1606.04076} multiplied by $y^{3/2}$ is computed as follows.

  \begin{lstlisting}
    sage: from homogeneous_space import *
    sage: from elliptic_genus import elliptic_genus
    sage: P = ParabolicSubgroup(CartanType('G2'), None, [1, 2])
    sage: X = HomogeneousSpace(P)
    sage: E = CompletelyReducibleEquivariantVectorBundle(X, [(2, 0), (0, 1), (0, 1)])
    sage: Y = CompleteIntersection(E)
    sage: elliptic_genus(Y, 3)
    -36*y - 36*y^2 + (36*y^-1 - 36*y - 36*y^2 + 36*y^4)*q + (36*y^-2 + 36*y^-1 - 72*y - 72*y^2 + 36*y^4 + 36*y^5)*q^2 + (36*y^-2 + 72*y^-1 - 108*y - 108*y^2 + 72*y^4 + 36*y^5)*q^3 + O(q^4)
  \end{lstlisting}

\end{example}

\section*{Acknowledgements}

The classes \lstinline|ParabolicSubgroup| and \lstinline|HomogeneousSpace| were co-developed with Akihito Nakamura and Kazushi Ueda.
I would like to express my gratitude to Kazushi Ueda,
my supervisor,
for raising issues about elliptic genera and suggesting to write codes in SageMath,
as well as for his thoughtful advice and encouragement during the writing of this paper.
I thank to Akihito Nakamura for his careful review of the docstring.
I'm also grateful to Makoto Miura for his meticulous review and bug detection.
My thanks go out to Kenichiro Kobayashi for the optimization idea regarding the calculation of the elliptic genus.

\bibliography{references}
\bibliographystyle{amsalpha}

\end{document}